\documentclass[12pt]{article}
\usepackage{amssymb}
\usepackage{amsmath,graphicx,amsfonts,epsfig,eucal}
\usepackage{amsfonts}

\newcommand{\ds}{\displaystyle}
\newcommand{\pa}{\partial}
\newcommand{\rr}{{\rm{I \! R}}}
\newtheorem{theorem}{Theorem}

\newtheorem{proposition}[theorem]{Proposition}

\input{tcilatex}

\begin{document}

\title{An economic game with stochastic dynamics}
\author{A. L. Ciurdariu$^{a}$, M. Neam\c{t}u$^{b}$, D. Opri\c s$^{c}$}
\date{}
\maketitle

\begin{tabular}{cccccccc}
\scriptsize{$^{a}$ Department of Mathematics,Politehnica University of Timi\c{s}oara}\\
\scriptsize{P-\c{t}a. Victoriei, nr 2,
300004, Timi\c{s}oara, Romania,}\\
\scriptsize{E-mail: cloredana43@yahoo.com}\\
\scriptsize{$^{b}$Department of Economic Informatics, Mathematics and Statistics,}\\
\scriptsize{Faculty of Economics, West University of Timi\c soara,}\\
\scriptsize{str. Pestalozzi, nr. 16A, 300115, Timi\c soara, Romania,}\\
\scriptsize{E-mail:mihaela.neamtu@fse.uvt.ro,}\\
\scriptsize{$^{c}$ Department of Applied
Mathematics, Faculty of Mathematics,}\\
\scriptsize{West University of Timi\c soara, Bd. V. Parvan, nr. 4, 300223, Timi\c soara, Romania,}\\
 \scriptsize{E-mail: opris@math.uvt.ro}\\

\end{tabular}

\begin{abstract} In this paper we investigate a
stochastic model for an economic game. To describe this model we
have used a Wiener process, as the noise has a stabilization effect.
The dynamics are studied in terms of stochastic stability in the
stationary state, by constructing the Lyapunov exponent, depending
on the parameters that describe the model. The numerical simulation
that we did justifies the theoretical results. \end{abstract}

\medskip




\section{Introduction.}

\qquad Stochastic modeling plays an important role in many branches
of science. In many practical situations, perturbations are
expressed in terms of white noise, modeled by brownian motion. The
behavior of a deterministic dynamical system which is disturbed by
noise may be modeled by a stochastic differential equation (SDE),
\cite{KP}. The stochastic stability has been introduced by Bertram
and Sarachik and is characterized by the negativeness of Lyapunov
exponents. In general, it is not possible to determine this
exponents explicitly. Many numerical approaches have been proposed,
which generally used the simulation of the stochastic trajectories
\cite{Jed}. In the present paper, we study a stochastic dynamical
system that are used in economy, in describing a Counot duopoly
game.

In 1838, Cournot introduced the first formal theory of oligopoly,
which treated the case of naive expectations, where each player
assumes the last values taken by the competitors without estimation
of their future reactions \cite{Cournot}. Recently, a lot of
articles have shown that the Cournot model may lead to a cyclic or
chaotic behavior \cite{Bundau}, \cite{Bischi}, \cite{Mircea},
\cite{Puu1}, \cite{Puu2}, \cite{Puu3}. Also, in \cite{Rosser},
Rosser reviews the development of the theory of complex oligopoly
dynamics.

In the present paper we have studied a stochastic Cournot economic
game. In Section 2 we present the Lyapunov exponent and stability in
stochastic 2d dynamical structures. Section 3 studies the Lyapunov
exponent for an economic game with stochastic dynamics. Some
numerical simulations are given in Section 4. Finally, Section 5
draws some conclusions.

\section{The Lyapunov exponent and stability in stochastic 2d dynamical structures.}

\qquad Let $(\Omega , {\cal F}, {\cal P})$ be a probability space
\cite{KP}. It is
 assumed that the $\sigma -$algebra ${\cal F}$ is a filtration that
 is, ${\cal F}$ is generated by a family of $\sigma -$algebra ${\cal F}_t(t\geq
 0)$ such that

$${\cal F}_s\subset {\cal F}_t\subset {\cal F}, \quad \forall s\leq t, s,t\in
I,$$ where $I=[0, T]$, $T\in (0, \infty)$.

Let $\{x(t)=(x_1(t), x_2(t))\}_{t\geq 0}$ be a stochastic process.
The system of Ito equations:
\begin{equation}\label{1}
dx_i(t,\omega)=f_i(t,x(t,\omega))dt+g_i(x(t,\omega))dw(t,\omega),
i=1,2,
\end{equation} with the initial condition $x(0)=x_0$ is written as:
\begin{equation}\label{2}
x_i(t,\omega)=x_{i0}(\omega)+\int_0^tf_i(x(s,\omega))ds+\int_0^tg_i(x(s,\omega))dw(s,\omega),
i=1,2,
\end{equation} for almost all $\omega\in\Omega $ and for each $t>0$,
where $f_i(x)$ are drift functions, $g_i(x)$ are diffusion
functions, $\int_0^tf_i(x(s))ds$, $i=1,2$ are Riemann integrals and
$\int_0^tg_i(x(s))dw(s)$ are It$\hat{o}$ integrals. It is assumed
that $f_i$ and $g_i$, $i=1,2$ satisfy the conditions of existence of
solution for this SDE with initial condition $x(0)=a_0\in\rr^n$.

Let $x_0=(x_{10}, x_{20})\in\rr^2$ be a solution of the system:
\begin{equation}\label{3}
f_i(x_0)=0, i=1,2.
\end{equation}

The functions $g_i, i=1,2$ are chosen so that:
\begin{equation*}
g_i(x_0)=0, i=1,2.
\end{equation*}

In what follows, we consider:
\begin{equation*}\label{4}
g_i(x)=\sum_{j=1}^{2}b_{ij}(x_j-x_{0j}), i=1,2,
\end{equation*} where $b_{ij}\in\rr, i,j=1,2.$

The linearized system of (\ref{2}) in $x_0$, is given by:
\begin{equation*}\label{5}
X(t)=\int_0^tAX(s)ds+\int_0^tBX(s)dw(s),
\end{equation*} where
\begin{equation*}\label{6}
X(t)=\left (\begin{array}{c} x(t,\omega )\\ y(t,\omega
)\end{array}\right ), A=\left (\begin{array}{cc} a_{11} & a_{12}\\
a_{21} & a_{22}\end{array}\right ), B=\left (\begin{array}{cc} b_{11} & b_{12}\\
b_{21} & b_{22}\end{array}\right ),
\end{equation*}
\begin{equation*}\label{7}
a_{ij}=\ds\frac{\pa f_i}{\pa x_j}|_{x_0}, b_{ij}=\ds\frac{\pa
g_i}{\pa x_j}|_{x_0}.
\end{equation*}

The Oseledec multiplicative ergodic theorem \cite{Ose} asserts the
existence of 2 non-random Lyapunov exponents
$\lambda_2\leq\lambda_1=\lambda$. The top Lyapunov exponent is given
by:
\begin{equation*}\label{8}
\lambda =\lim_{t\to\infty} \sup\log\sqrt{x(t)^2+y(t)^2}.
\end{equation*}

Applying the change to polar coordinates:
\begin{equation*}
x(t)=r(t)cos \theta (t), y(t)=r(t) sin\theta (t)
\end{equation*} by writing the It$\hat{o}$ formula for
\begin{equation*}
h_1(x,y)=\ds\frac{1}{2}\log (x^2+y^2)=\log (r), h_2(x,y)=arctg
(\ds\frac{y}{x})=\theta .
\end{equation*} we get:

\begin{proposition}\label{P1} \cite{Jed}.  The formulas
\begin{equation}\label{11}
\log \left (\ds\frac{r(t)}{r(0)}\right
)\!=\!\int_0^tq_1(\theta(s))\!+\!\ds\frac{1}{2}(q_4(\theta(s))^2\!-\!q_2(\theta(s))^2)ds\!+\!\int_0^tq_2(\theta(s))dw(s),
\end{equation}
\begin{equation}\label{12}
\theta(t)\!=\!\theta(0)+\int_0^tq_3(\theta(s))\!-\!q_2(\theta(s)q_4(\theta(s))ds\!+\!\int_0^tq_4(\theta(s))dw(s),
\end{equation} hold,
 where
\begin{equation}\label{13}
\begin{array}{llll}
q_1(\theta)=a_{11}cos^2(\theta)+(a_{12}+a_{21})cos\theta\sin\theta+a_{22}sin^2\theta,\\
q_2(\theta)=b_{11}cos^2(\theta)+(b_{12}+b_{21})cos\theta\sin\theta+b_{22}sin^2\theta,\\
q_3(\theta)=a_{21}cos^2(\theta)+(a_{22}-a_{11})cos\theta\sin\theta-a_{12}sin^2\theta,\\
q_4(\theta)=b_{21}cos^2(\theta)+(b_{22}-b_{11})cos\theta\sin\theta-b_{12}sin^2\theta.\\
\end{array}
\end{equation}
\end{proposition}
As the expectation of the It$\hat{o}$ stochastic integral is null
$$E\int_0^tq_2(\theta (s))dw(s)=0,$$ the Lyapunov exponent is given
by:
$$\lambda\!=\!\lim_{t\to\infty}\ds\frac{1}{t}\log\left (\ds\frac{r(t)}{r(0)}\right )\!=\!
\lim_{t\to\infty}\ds\frac{1}{t}E\int_0^t(q_1(\theta(s))\!+\!\ds\frac{1}{2}(q_4(\theta(s)))^2\!-\!q_2(\theta(s)))ds.$$
Applying the Oseledec theorem, if $r(t)$ is ergodic, we get:
\begin{equation*}\label{14}
\lambda
=\int_0^t(q_1(\theta)+\ds\frac{1}{2}(q_3(\theta)^2-q_2(\theta)))p(\theta)d\theta
,
\end{equation*}where $p(\theta)$ is the probability distribution of
the process $\theta$.

An approximation of this distribution is calculated by solving the
Fokker-Planck equation.

The Fokker-Planck (FPE) equation associated with equation (\ref{12})
for $p=p(t,\theta)$ is
\begin{equation}\label{15}
\ds\frac{\pa p}{\pa t}+\ds\frac{\pa }{\pa\theta
}((q_3(\theta)-q_2(\theta)q_4(\theta))p)-\ds\frac{1}{2}\ds\frac{\pa
^2}{\pa \theta^2}(q_4(\theta)^2p)=0.
\end{equation}

From (\ref{15}), it results that the solution $p(\theta)$ of the FPE
is solution of the following first order equation:
\begin{equation}\label{16}
(-q_3(\theta)\!+\!q_1(\theta)q_4(\theta)\!+
\!q_2(\theta)g_5(\theta))p(\theta)\!+\!\ds\frac{1}{2}q_4(\theta)^2p'(\theta)\!=\!p_0,
\end{equation} where $p'(\theta )=\ds\frac{dp}{d\theta}$ and
\begin{equation*}\label{17}
q_5(\theta)=-(b_{12}+b_{21})sin 2\theta -(b_{22}-b_{11})cos 2\theta
.
\end{equation*}

\begin{proposition}\label{prop2}\cite{Jed}. If $q_4(\theta)\neq 0$, the
solution of the equation (\ref{16}) is given by:
\begin{equation*}\label{18}
p(\theta)=\ds\frac{k}{D(\theta)q_4(\theta)^2}\left (1+\eta
\int_0^\theta D(u)du\right )
\end{equation*} where $k$ is determined by the normality condition
\begin{equation*}\label{18}
\int_0^{2\pi }p(\theta )d\theta =1
\end{equation*} and
\begin{equation*}\label{19}
\eta =\ds\frac{D(2\pi )-1}{\int_0^{2\pi }D(u)du}.
\end{equation*}
The function $D$ is given by:
\begin{equation*}\label{21}
D(\theta )=\exp
(-2\int_0^\theta\ds\frac{q_3(u)-q_2(u)q_4(u)-q_4(u)q_5(u)}{q_4(u)^2}du)
\end{equation*}
\end{proposition}

A numerical solution of the phase distribution could be performed by
a simple backward difference scheme.

We consider $N\in\rr_+$, $h=\ds\frac{\pi}{N}$ and
\begin{equation*}\label{13}
\begin{array}{llll}
q_1(i)=a_{11}\cos^2(ih)+(a_{12}+a_{21})\cos(ih)\sin(ih)+a_{22}\sin^2(ih),\\
q_2(i)=b_{11}\cos^2(ih)+(b_{12}+b_{21})\cos(ih)\sin(ih)+b_{22}\sin^2(ih),\\
q_3(i)=a_{21}\cos^2(ih)+(a_{22}-a_{11})\cos(ih)\sin(ih)-a_{12}sin^2(ih),\\
q_4(i)=b_{21}\cos^2(ih)+(b_{22}-b_{11})\cos(ih)\sin(ih)-b_{12}\sin^2(ih),\\
q_5(i)=-(b_{12}+b_{21})\sin(2ih)-(b_{22}-b_{11})cos(2ih), i=0,...,N\\
\end{array}
\end{equation*}

The function $p(i), i=0,...,N$ is given by the following relations:
$$p(i)=(p(0)+\ds\frac{q_4(i)^2p(i-1)}{2h})F(i)$$ where
$$F(i)=\ds\frac{2h}{2h(-q_3(i)+q_2(i)q_4(i)+q_4(i)q_5(i))+q_4(i)^2}.$$
The Lyapunov exponent is $\lambda =\lambda (N)$, where
$$\lambda (N)=\sum_{i=0}^{N}(q_1(i)+\ds\frac{1}{2}(q_4(i)^2-q_2(i)^2))p(i)h.$$

\begin{proposition}\label{prop3} If the matrix B is given by:
$$b_{11}=\alpha, b_{12}=-\beta, b_{21}=\beta, b_{22}=\alpha $$
then
\begin{equation*}
\begin{split}
p(\theta)&=\ds\frac{k}{\beta^2}\exp\{\ds\frac{1}{\beta^2}((a_{21}-a_{12}-\alpha\beta)\theta+\ds\frac{1}{2}(a_{11}-a_{22})\cos
2\theta+ \ds\frac{1}{2}(a_{21}-\\
&-a_{12})\sin2\theta)\}
\end{split}
\end{equation*}
\begin{equation*}
k\!=\!\ds\frac{\beta^2}{\int_0^{2\pi}\exp\{\ds\frac{1}{\beta^2}((a_{21}\!\!-\!\!a_{12}\!\!-\!\!\alpha\beta)\theta\!+\!\ds\frac{1}{2}(a_{11}\!-\!a_{22})\cos
2\theta\!\!+\!\!
\ds\frac{1}{2}(a_{21}\!\!-\!\!a_{12})\sin2\theta)d\theta}
\end{equation*}
\begin{equation*}
\lambda=\ds\frac{1}{2}(a_{11}+a_{22}+\beta^2-\alpha^2)+\ds\frac{1}{2}(a_{11}-a_{22})c_2+\ds\frac{1}{2}(a_{21}+a_{12})s_2,
\end{equation*} where
$$c_2=\int_0^{2\pi}cos(2\theta)p(\theta)d\theta, \quad s_2=\int_0^{2\pi} sin(2\theta)p(\theta)d\theta .$$
\end{proposition}

\section{The Lyapunov exponent for an economic game with stochastic dynamics.}
\qquad Two firms enter the market with a homogenous consumption
product. The elements which describe the model are: the quantities
which enter the market
from the two firms $x_{i}\geq 0,$ $i=\overline{1,2};$  the inverse demand function $p:\mathbb{R}%
_{+}\rightarrow \mathbb{R}_{+}$ ( $p$ is a derivable function with $%
p^{\prime }\left( x\right) <0,\underset{x\rightarrow
a_{1}}{lim}p\left( x\right) =0,$ $\underset{x\rightarrow
0}{lim}p\left( x\right) =b_{1},\left( a_{1}\in
\overline{\mathbb{R}},b_{1}\in \overline{\mathbb{R}}\right) $;  the cost functions $C_{i}:%
\mathbb{R}_{+}\rightarrow \mathbb{R}_{+}$ ( $C_{i}$ are derivable
functions with $C_{i}^{\prime }\left( x_{i}\right) >0,$
$C_{i}^{\prime \prime }\geq 0,$ $i=\overline{1,2}$ ).

In our study we consider $p(x)=\ds\frac{1}{x}, x>0$ and
$C_i(x_i)=c_ix_i+d_i, i=1,2.$

The mathematical model of the stochastic dynamic economic game is
described by the stochastic system of equations:
\begin{equation}\label{31}
\begin{split}
x_1(t)\!=\!x_1\!(0)\!\!+\!k_1\!\int_0^t\!(\ds\frac{x_2(s)}{(x_1(s)\!+\!x_2(s))^2}\!-\!c_1)ds\!+\!\int_0^t(b_{11}x_1(s)\!+\!b_{12}x_2(s)\!+\!\gamma_1)dw(s)\\
x_2(t)\!=\!x_2\!(0)\!\!+\!k_2\!\int_0^t\!(\ds\frac{x_1(s)}{(x_1(s)\!+\!x_2(s))^2}\!-\!c_2)ds\!+\!\int_0^t(b_{21}x_1(s)\!+\!b_{22}x_2(s)\!+\!\gamma_2)dw(s)
\end{split}
\end{equation} where $b_{ij}\in\rr$, $i,j=1,2$, $k_1>0, k_2>0$, $x_i(t)=x_i(t,\omega
)$, $i=1,2$.
\begin{equation*}\label{32}
\gamma_1=-\ds\frac{b_{11}c_2+b_{12}c_1}{(c_1+c_2)^2},
\gamma_2=-\ds\frac{b_{21}c_2+b_{22}c_1}{(c_1+c_2)^2}.
\end{equation*}

For $b_{ij}=0$, $i,j=1,2$ model (\ref{31}) is reduced to the
classical model of the economic game \cite{Bundau}, \cite{Mircea}.

The system of stochastic equations (\ref{31}), (SDE), has the form
(\ref{2}) from section 2, where:
\begin{equation*}\label{33}
\begin{split}
f_1(x_1,x_2)=\ds\frac{x_2}{(x_1+x_2)^2}-c_1,
g_1(x_1,x_2)=b_{11}x_1+b_{12}x_2+\gamma_1,\\
f_2(x_1,x_2)=\ds\frac{x_1}{(x_1+x_2)^2}-c_2,
g_2(x_1,x_2)=b_{21}x_1+b_{22}x_2+\gamma_2.
\end{split}
\end{equation*}

Applying the results from section 2, we have:

\begin{proposition}\label{prop4}
(i) The stationary state of (SDE)
(\ref{31}) is given by:
\begin{equation*}\label{34}
x_{10}=\ds\frac{c_2}{(c_1+c_2)^2},
x_{20}=\ds\frac{c_1}{(c_1+c_2)^2};
\end{equation*}

(ii) The elements of the matrix $A$, which characterize linearized
equation (\ref{31}) in $(x_{10}, x_{20})$ are:
\begin{equation*}\label{35}
\begin{split}
a_{11}=-2k_1c_1(c_1+c_2), a_{12}=-k_1(c_1^2-c_2^2)\\
a_{21}=k_2(c_1^2-c_2^2), a_{22}=-2k_2c_2(c_1+c_2);
\end{split}
\end{equation*}

(iii) The roots of the characteristic equation:
\begin{equation}\label{36}
\mu^2-(a_{11}+a_{22})\mu+a_{11}a_{22}-a_{12}a_{21}=0
\end{equation} have the real part:
\begin{equation*}\label{37}
Re(\mu_{1,2})=-(k_1c_1+k_2c_2)(c_1+c_2);
\end{equation*}

(iv) If $b_{11}=\alpha$, $b_{12}=-\beta$, $b_{21}=\beta$,
$b_{22}=\alpha$, $\beta\neq 0$, then the Lyapunov coefficient of
(SDE) (\ref{3}) is:
\begin{equation}\label{38}
\begin{split}
\lambda\!
&=\!-(k_1c_1+k_2c_2)(c_1+c_2)+\ds\frac{1}{2}(\beta^2-\alpha^2)-
(k_1c_1-k_2c_2)(c_1+c_2)D_2+\\
&+\ds\frac{1}{2}(k_2-k_1)(c_1^2-c_2^2)E_2
\end{split}
\end{equation} where
\begin{equation*}\label{39}
D_2=\int_0^{2\pi}cos(2\theta)p(\theta)d\theta,
E_2=\int_0^{2\pi}sin(2\theta)p(\theta)d\theta
\end{equation*} and
\begin{equation*}\label{40}
\begin{split}
& p(\theta)=kg(\theta),
k=\ds\frac{1}{\int_0^{2\pi}g(\theta)d\theta},\\
&
g(\theta)=\ds\frac{1}{\beta^2}\exp\{\ds\frac{1}{\beta^2}((k_1+k_2)(c_1^2-c_2^2)+\alpha\beta)\theta-
(k_1c_1-k_2c_2)(c_1+c_2)\cos(2\theta)+\\
&+\ds\frac{1}{2}(k_1+k_2)(c_1^2-c_2^2)\sin(2\theta)\}.
\end{split}
\end{equation*}

\end{proposition}

\section{Numerical Simulations.}

\qquad We have done the numerical simulations using a program in
Maple 12. For $c_1=0.2$, $c_2=2$, $k_1=0.2$, $k_2=0.4$, $\beta=2$,
in figure 1 is displayed $(\alpha, \lambda (\alpha))$, where
$\lambda(\alpha)$ is given by (\ref{38}). For $\alpha\in (-\infty,
-1.2) \cup (1.1, \infty)$, the Lyapunov exponent is negative, then
(SDE) has an asymptotically stable stationary state. For $\alpha\in
(-1.2, 1.1)$, the Lyapunov exponent is positive and (SDE) has an
asymptotically unstable stationary state.

\begin{center}\begin{tabular}{ccc}
\\ Fig 1. $(\alpha, \lambda(\alpha))$\\
\includegraphics[width=5cm]{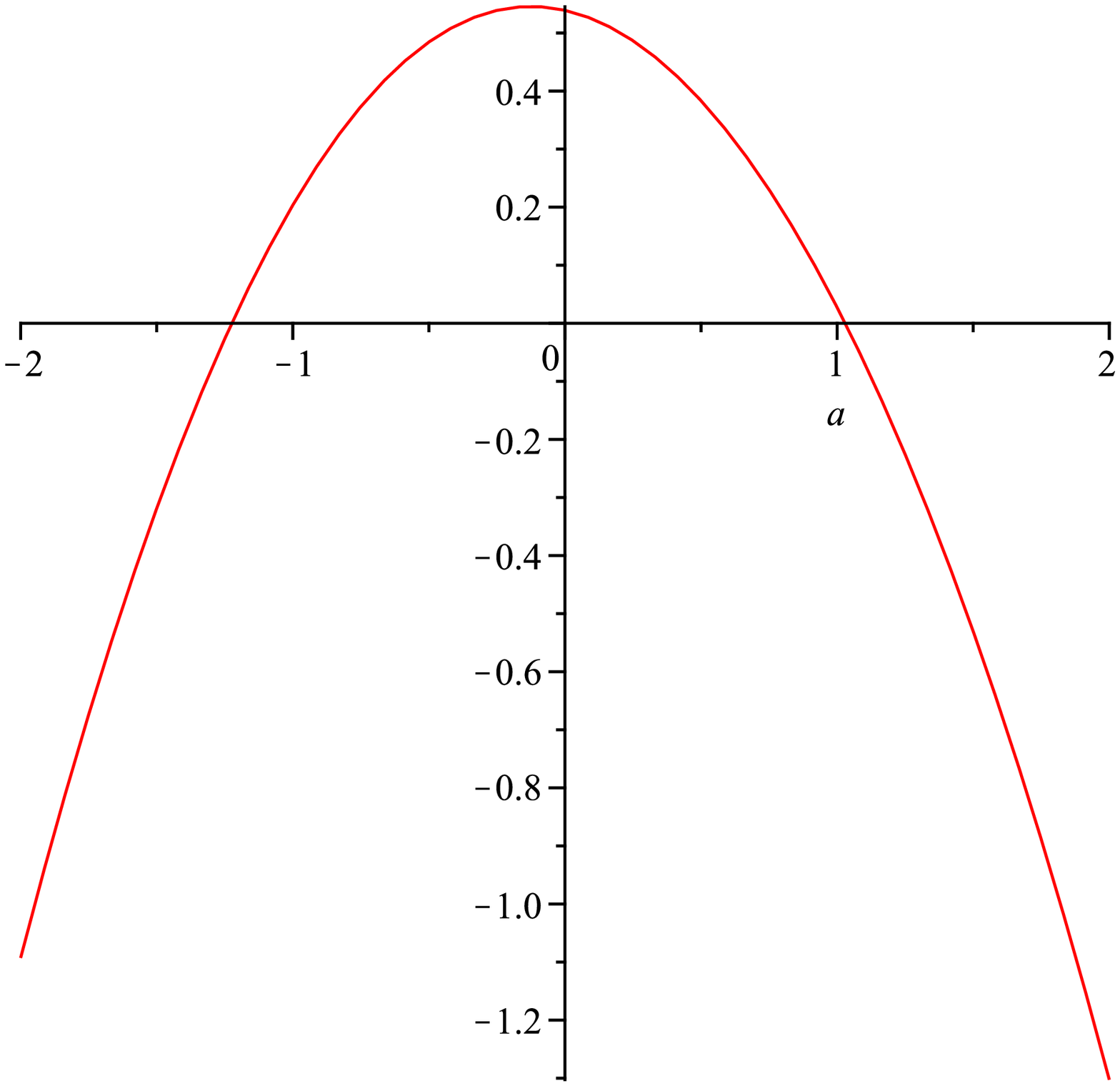}

\end{tabular}
\end{center}

If $\beta$ is a real parameter and $\alpha=2$, the figure 2 shows
the behavior of the top Lyapunov exponent as a function of $\beta$:
$(\beta, \lambda(\beta))$.

\begin{center}\begin{tabular}{ccc}
\\ Fig 2. $(\beta, \lambda(\beta))$\\
\includegraphics[width=5cm]{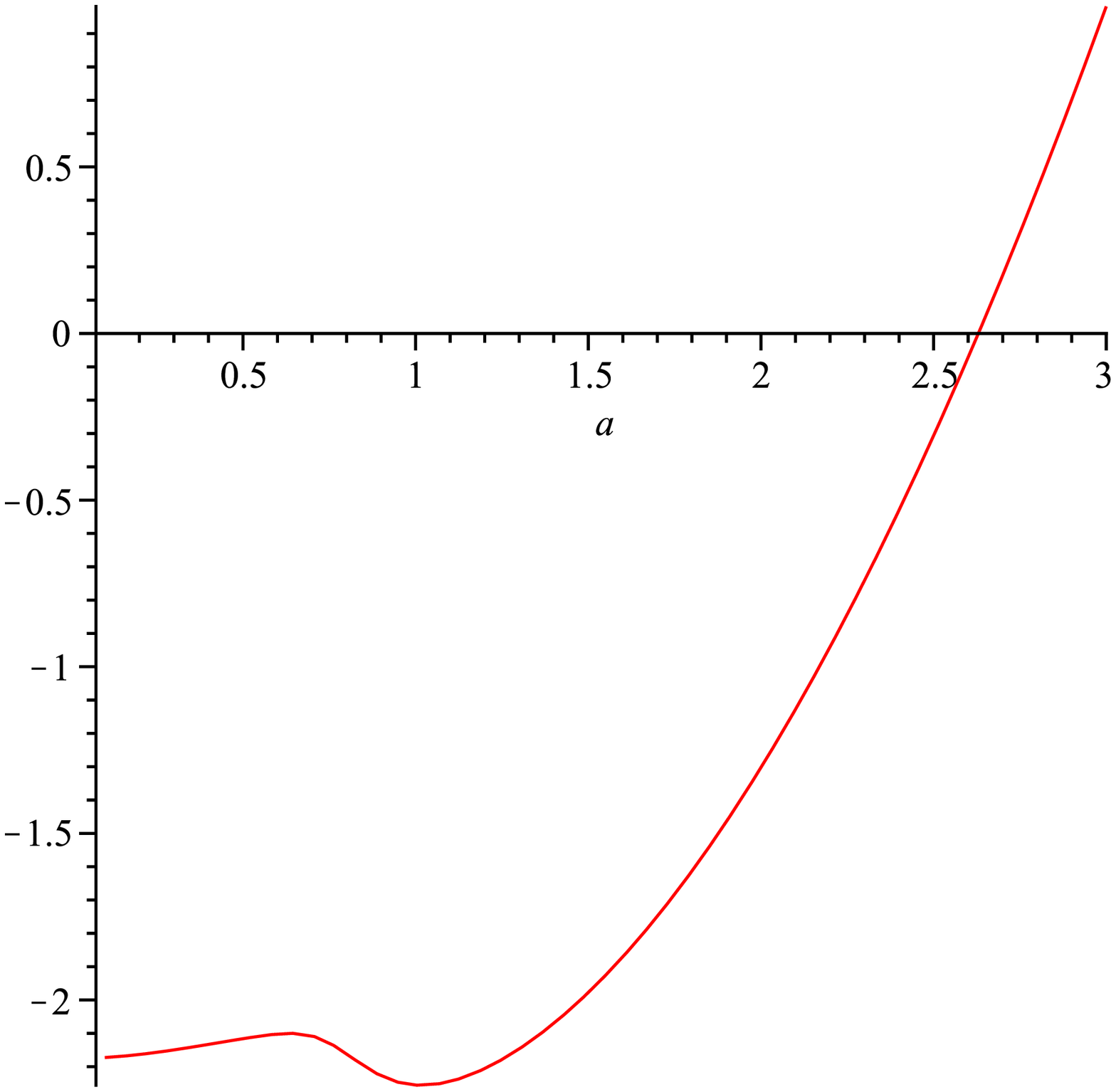}

\end{tabular}
\end{center}

For $\beta \in (-\infty, -2.6)\cup (2.6,\infty)$ the Lyapunov
exponent is positive and (SDE) has an asymptotically unstable
stationary state. For $\beta (-2.6, 2.6)$ the Lyapunov exponent is
negative and (SDE) has an asymptotically stable stationary state.

The Euler second order scheme for (SDE) (2) is given by:
\begin{equation*}
\begin{split}
&x_1(n+1)=x_1(n)+h\left (\ds\frac{x_2(n)}{(x_1(n)+x_2(n))^2}-c_1\right )+(b_{11}x_1(n)+b_{12}x_2(n)+\gamma_1)\\
&\cdot
G(n)+b_{11}(b_{11}x_1(n)+b_{12}x_2(n)+\gamma_1)\ds\frac{G(n)^2-h}{2}+
(-\ds\frac{2x_1(n)x_2(n)}{(x_1(n)+x_2(n))^3}\\
&\cdot\left (\ds\frac{x_2(n)}{(x_1(n)+x_2(n))^2}-c_1\right
)+(b_{11}x_1(n)+b_{12}x_2(n)+\gamma_1)\ds\frac{x_1(n)x_2(n)}{(x_1(n)+x_2(n))^3}
)\ds\frac{h^2}{2}+\\
&(b_{11}-\ds\frac{2x_2(n)}{(x_1(n)+x_2(n))^3})(b_{11}x_1(n)+b_{12}x_2(n)+\gamma_1)\ds\frac{hG(n)}{2},
\end{split}
\end{equation*}

\begin{equation*}
\begin{split}
&x_2(n+1)=x_2(n)+h\left (\ds\frac{x_1(n)}{(x_1(n)+x_2(n))^2}-c_2\right )+(b_{21}x_1(n)+b_{22}x_2(n)+\gamma_2)\\
&\cdot
G(n)+b_{22}(b_{21}x_1(n)+b_{22}x_2(n)+\gamma_2)\ds\frac{G(n)^2-h}{2}+
(-\ds\frac{2x_1(n)x_2(n)}{(x_1(n)+x_2(n))^3}\\
&\cdot\left (\ds\frac{x_1(n)}{(x_1(n)+x_2(n))^2}-c_2\right
)+(b_{21}x_1(n)+b_{22}x_2(n)+\gamma_2)\ds\frac{x_1(n)x_2(n)}{(x_1(n)+x_2(n))^3}
)\ds\frac{h^2}{2}+\\
&(b_{21}-\ds\frac{2x_1(n)}{(x_1(n)+x_2(n))^3})(b_{21}x_1(n)+b_{22}x_2(n)+\gamma_2)\ds\frac{hG(n)}{2},
\end{split}
\end{equation*} where $G(n)=w((n+1)h)-w(nh)$, $n=1,2,...$, and $x_i(n)=x_i(nh,\omega)$,
 $i=1,2$.

In figures 3 and 4 are displayed the orbits: $(n, x_1(n, \omega))$
for (SDE) and $(n, x_1(n))$ for (ODE):

\begin{center}\begin{tabular}{ccc}
\\ Fig 3. $(n, x_1(n, \omega))$ & Fig 4. $(n, x_1(n))$\\
\includegraphics[width=5cm]{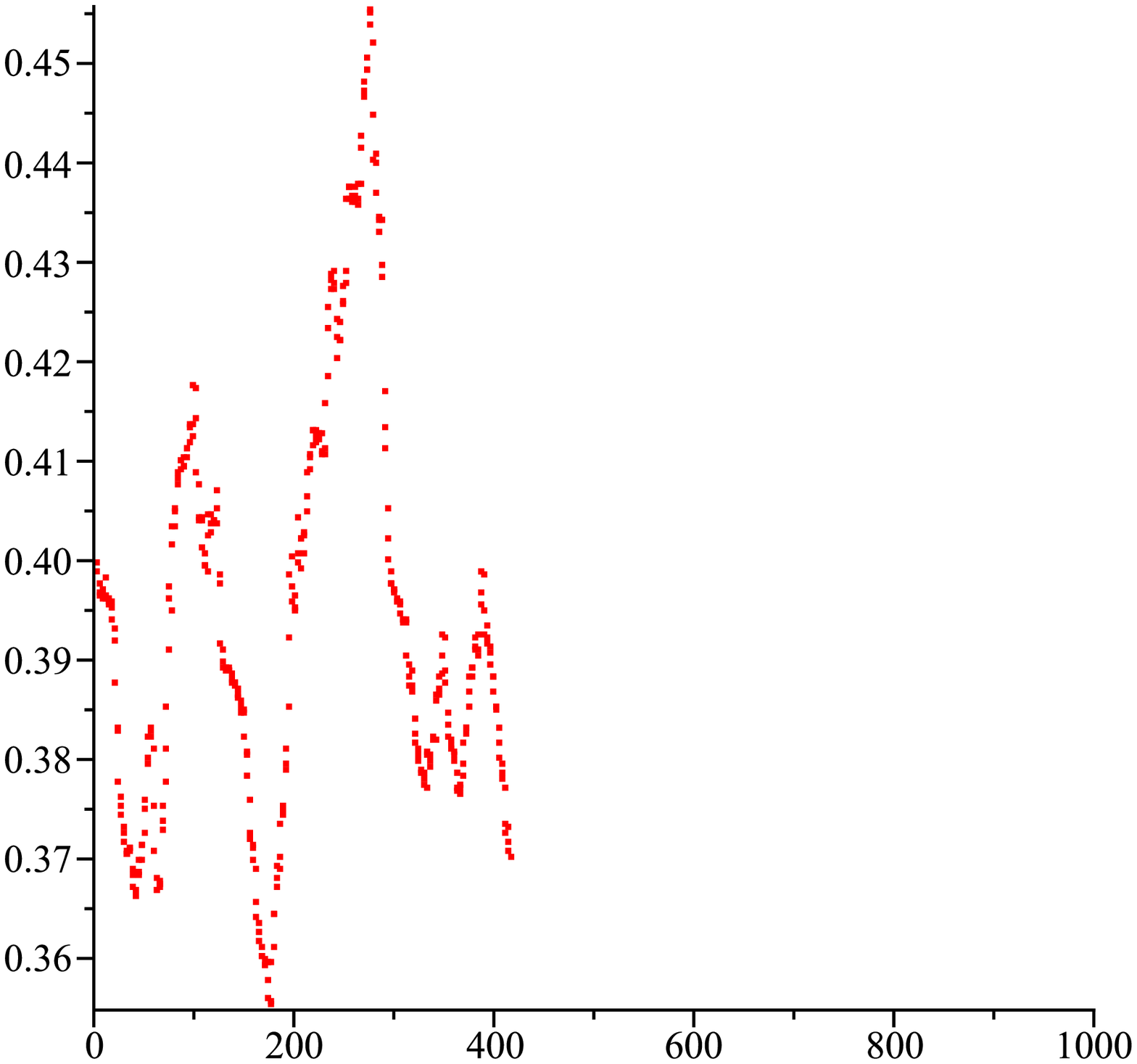} & \includegraphics[width=5cm]{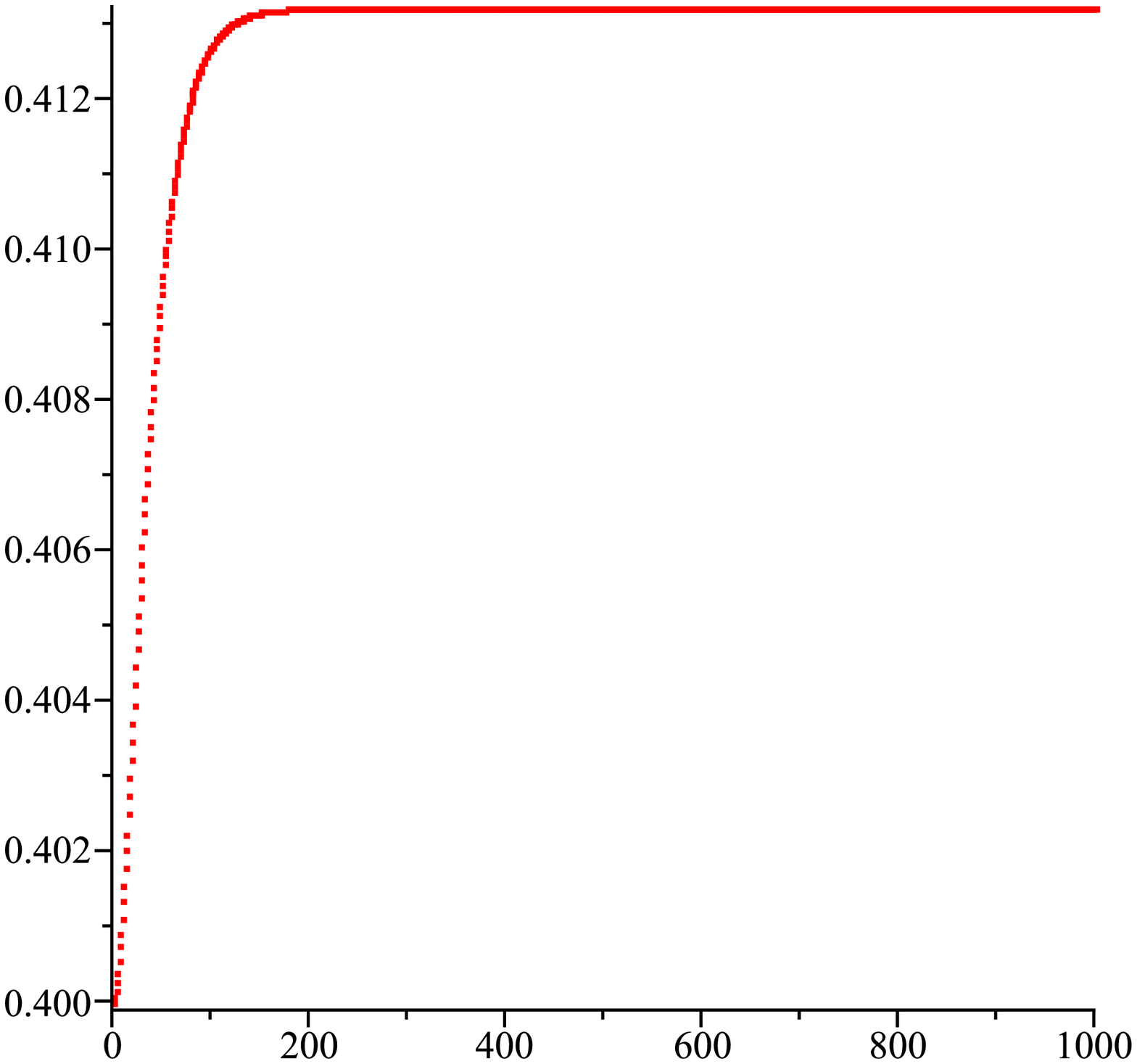}

\end{tabular}
\end{center}

In figures 5 and 6 are displayed the orbits: $(n, x_2(n, \omega))$
for (SDE) and $(n, x_2(n))$ for (ODE):

\begin{center}\begin{tabular}{ccc}
\\ Fig 5. $(n, x_2(n, \omega))$ & Fig 6. $(n, x_2(n))$\\
\includegraphics[width=5cm]{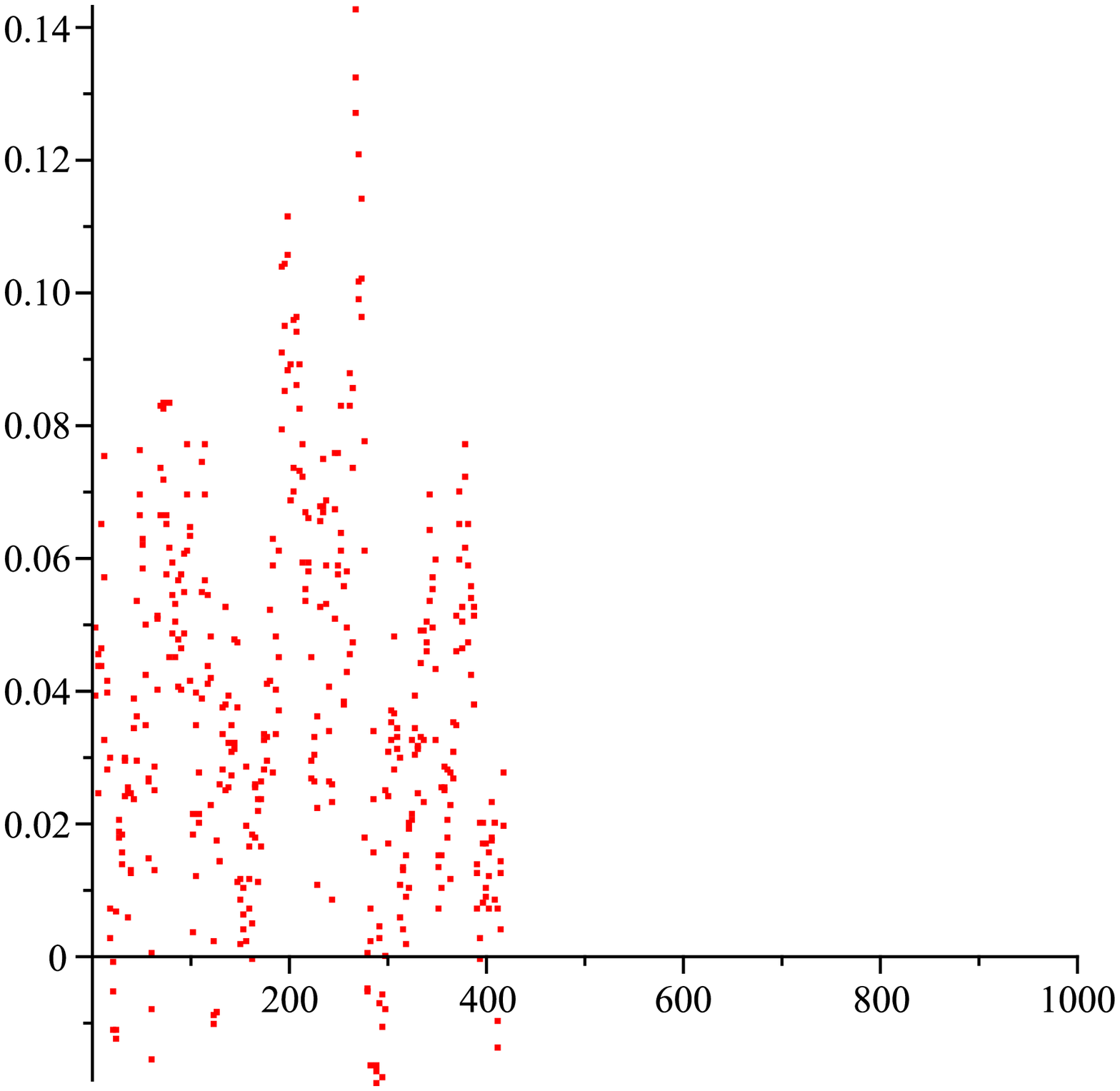} & \includegraphics[width=5cm]{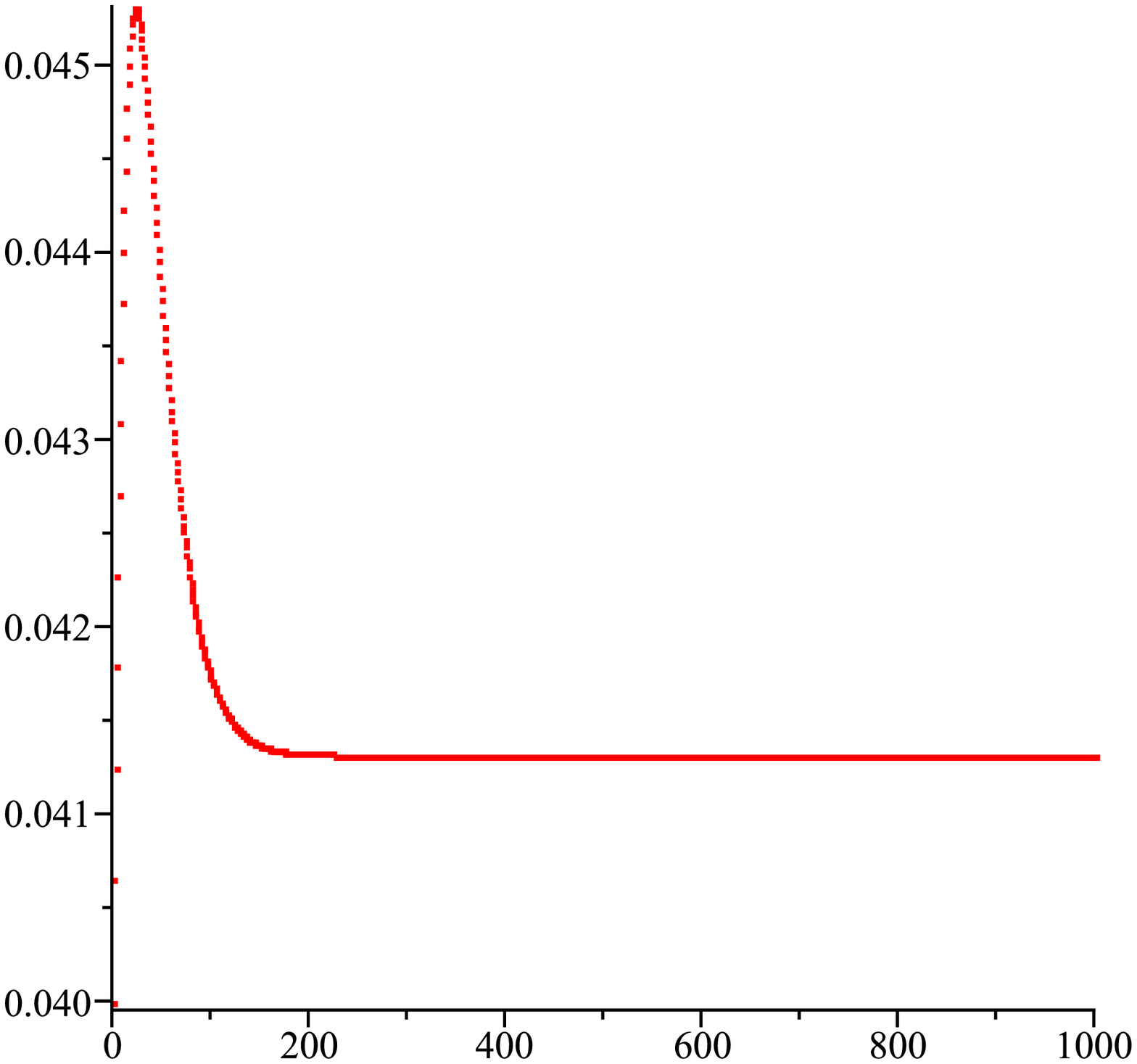}

\end{tabular}
\end{center}

In figures 7 and 8 are displayed the orbits: $(x_1(n,\omega), x_2(n,
\omega))$ for (SDE) and  $(x_1(n), x_2(n))$ for (ODE):

\begin{center}\begin{tabular}{ccc}
\\ Fig 7. $(x_1(n,\omega), x_2(n,
\omega))$ & Fig 8. $(x_1(n), x_2(n))$\\
\includegraphics[width=5cm]{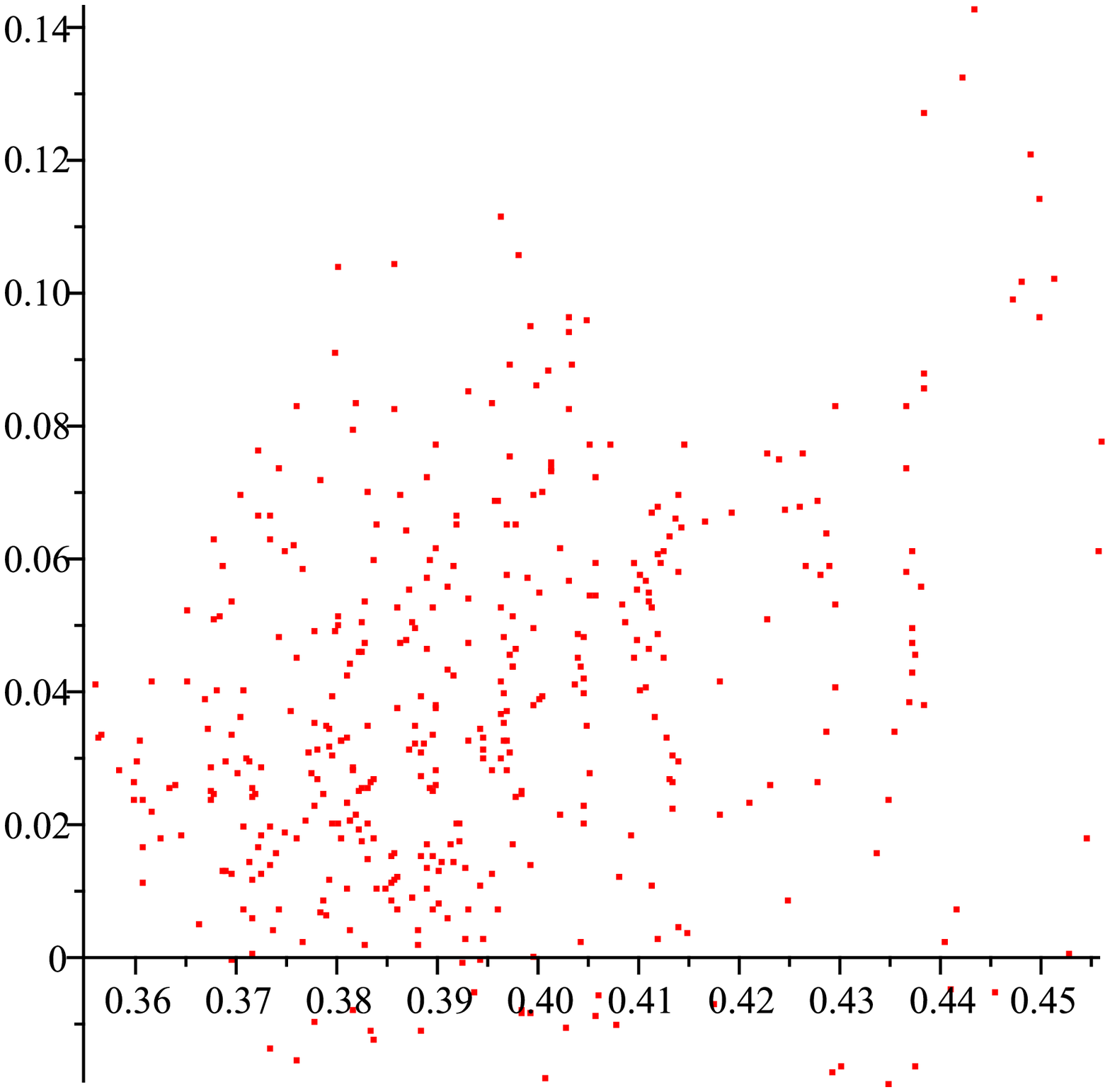} & \includegraphics[width=5cm]{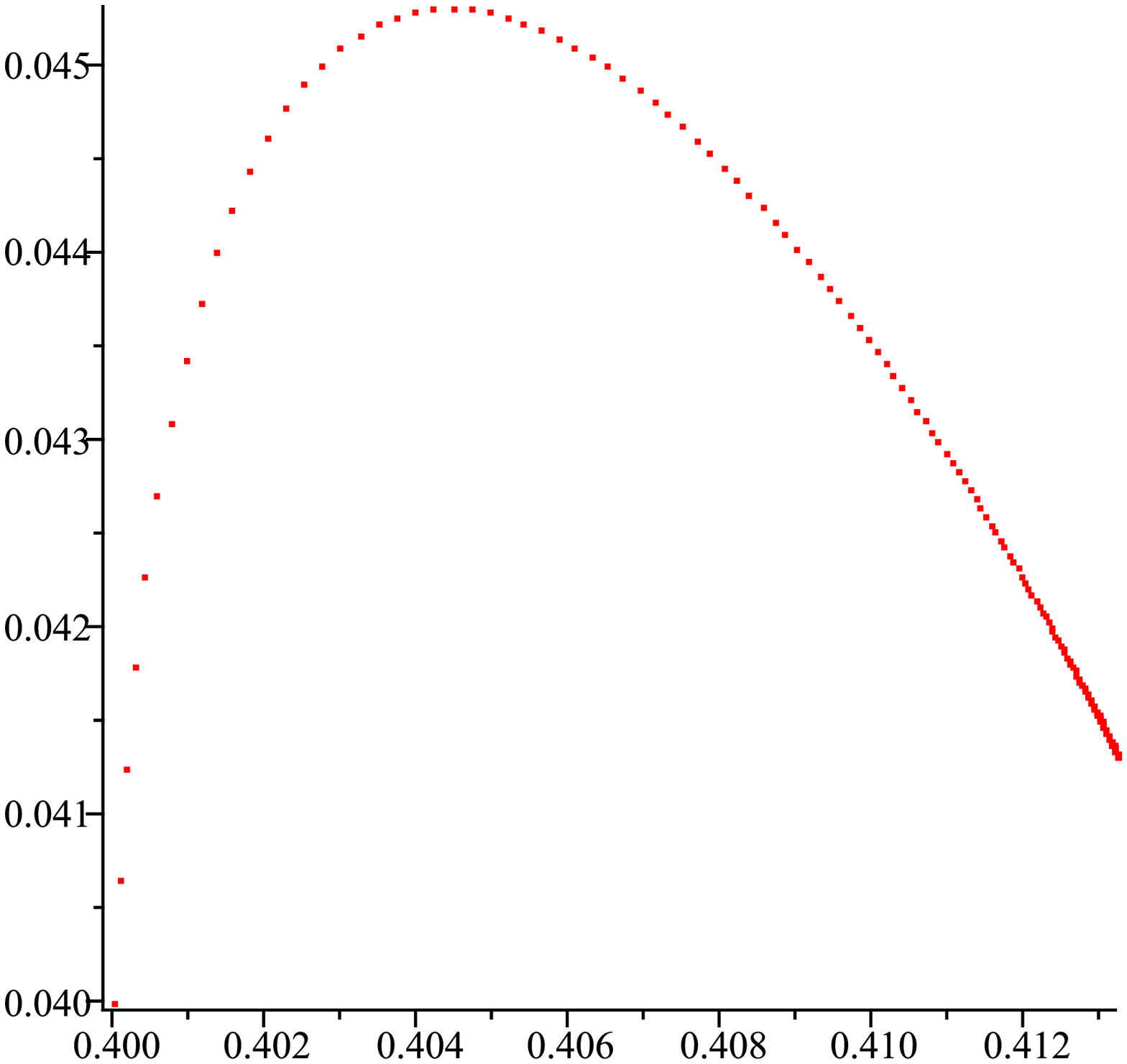}

\end{tabular}
\end{center}

\section{Conclusions.}

\qquad In the present paper we investigate an economic game with
stochastic dynamics. We focus on a particular game and determine the
Lyapunov exponent for the stochastic system of equations that
describes the mathematical model. The calculation of the top
Lyapunov exponent enables us to decide whether a stochastic system
is stable or not. Using a program in Maple 12, we display the
Lyapunov exponent and the system orbits.

\end{document}